\documentclass[english, 12pt]{article}
\usepackage[english]{babel,varioref}
\usepackage[ansinew]{inputenc}
\usepackage[fleqn,reqno]{amsmath}
\usepackage{paralist}
\usepackage{natbib}
\usepackage{graphics}
\usepackage{multirow}
\usepackage{setspace}

\bibpunct[: ]{(}{)}{,}{a}{}{,}

 \setdefaultenum{(i)}{a)}{(1)}{(A)}

\newtheorem{thm}{Theorem}
\newtheorem{cor}{Corollary}

\def\max{\mathop{\rm max}}

\begin{document}

\title{Comparing Different Information Levels}

\author{Uwe Saint-Mont\footnote{
Nordhausen University of Applied Sciences, Germany}}

\maketitle



{\bf AMS-Classification:} 94A17, 93C41, 60G40

\vspace{10ex}

{\bf Abstract}. Given a sequence of random variables ${\bf X}=X_1,X_2,\ldots$
suppose the aim is to maximize one's return by picking a `favorable' $X_i$.
Obviously, the expected payoff crucially depends on the information at hand. An
optimally informed person knows all the values $X_i=x_i$ and thus receives $E (\sup
X_i)$. We will compare this return to the expected payoffs of a number of observers
having less information, in particular $\sup_i (EX_i)$, the value of the sequence to
a person who only knows the first moments of the random variables.

In general, there is a stochastic environment (i.e. a class of random variables $\cal C$), and several levels of information. Given some ${\bf X} \in {\cal C}$, an observer possessing information $j$ obtains $r_j({\bf X})$. We are going to study `information sets' of the form
$$
R_{\cal C}^{j,k} = \{ (x,y) | x = r_j({\bf X}), y=r_k({\bf X}), {\bf X} \in {\cal C}
\},
$$
characterizing the advantage of $k$ relative to $j$. Since such a set measures
the additional payoff by virtue of increased information, its analysis yields a number of interesting results, in particular
`prophet-type' inequalities.

\newpage

\section{Several Information Levels}

Suppose there is a sequence of bounded random variables ${\bf X} =X_1,X_2,\ldots$
and the aim is to maximize one's return by picking a `favorable' $X_i$.
The first aim of this contribution is to study observers with different kinds of
information:

Suppose an observer knows all the realizations of the random variables and may thus
choose the largest one. His expected return is therefore
\begin{equation}\label{prophet}
m= M({\bf X})=  M(X_1,X_2,\ldots) =  E(\sup_i X_i),
\end{equation}
which is called the value to a {\it prophet}. Since the prophet always picks the
largest realization his value $m$ is a natural upper bound, given a sequence $\bf
X$.

Traditionally, $m$ has been compared to the value obtained by a {\it statistician} who
observes the process sequentially. This gambler, studied in detail in \citet{ch71},
relies on stopping rules $T$, which have to be measurable with respect to the
$\sigma$-field of past events. Behaving optimally the statistician may thus receive
\begin{equation}\label{ext}
v=V({\bf X}) = \sup_{T} EX_T,
\end{equation}
If there is a finite horizon $n$, one defines $v=\sup_{T \in {\cal T}, T \le n}
EX_T$ and $m=E(\max_{1\le i \le n} X_i).$ To avoid trivialities, we assume $n \ge 2$
throughout this article.

A {\it minimally informed gambler} has to make his choice on the basis that he only
knows the random variables' expected values. Behaving optimally, he gets
\begin{equation}
u=U({\bf X})=\sup_i EX_i,
\end{equation}
an amount that is entirely due to his (weak) prior information, and is a
straightforward counterpart to $E (\sup_i X_i)$.

One might think that a person who knows the {\it common distribution}
$L(X_1,\ldots,X_n)$ (but none of the observations) should receive a larger payoff.
However, no matter how this gambler makes up his mind, at the end of the day he has
to choose an index $i \in \{1,\ldots,n \}$, and thus his expected reward will be
largest if $EX_i=u$. Thus, although he knows much more than the minimally informed
gambler his superior knowledge does not pay off.

In other words, it's the observations that make a difference. Suppose a person knows
the dependence structure among the random variables {\it and} some of the
observations, w.l.o.g. $x_1,\ldots,x_j$. Notice, that there is no sequential
unfolding of information, however, this {\it partially informed gambler} may use the
values known to him to update his knowledge on the variables not observed, i.e. he
may refer to conditional expectations. Thus he obtains
$$
\max [ x_1,\ldots,x_j,E(X_{j+1} | x_1,\ldots,x_j),\ldots,E(X_{n} | x_1,\ldots,x_j)],
$$
and his expected return is
\begin{equation}
w=W({\bf X})=E ( \max(X_1,\ldots,X_j,E(X_{j+1} | X_1,\ldots,X_j),\ldots,E(X_{n} |
X_1,\ldots,X_j)) ).
\end{equation}

This observer can be reduced to a classical situation as follows: Given
$x_1,\ldots,x_j$, he will only consider the largest of these values; and the same
with $E(X_{j+1}| x_1,\ldots,x_j),$ $\ldots,$ $E(X_{n}| x_1,\ldots,x_j)$. Thus
w.l.o.g. it suffices to compare $$\max(x_1,\ldots,x_j) \;\;\; \mbox{and} \;\;\;
\max(E(X_{j+1}| x_1,\ldots,x_j),\ldots, E(X_{n}| x_1,\ldots,x_j))$$ which is
tantamount to the comparison of $v$ and $m$ if $n=2$ and if arbitrary dependencies
are allowed. In this situation the statistician behaves optimally if he chooses
$x_1$ whenever $x_1 \ge E(X_2|x_1)$. Thus, the set of all possible values here is
given by $\{(w,m) | w \le m \le 2w-w^2, 0 \le w \le 1\}$ if w.l.o.g. $0 \le X_i \le
1$ for all $i$.


\subsection*{Stochastic environments (classes of random variables)}

For some fixed $\bf X$, the difference between two observers with different amounts
of information can be nonexistent or arbitrarily large. In order to quantify the
``value'' of information it is thus necessary to shift attention to some class of
random variables $\cal C$, where $M({\bf X})$ is finite (and nonnegative) for all
${\bf X} \in {\cal C}$. It is then natural to consider the worst case scenarios.
Traditionally these have been called prophet inequalities $ M({\bf X})-V({\bf X})
\le a$ and $M({\bf X})/V({\bf X}) \le b $ with smallest possible constants $a$ and
$b$ that hold for all ${\bf X} \in {\cal C}$.

Such stochastic inequalities follow easily from the more fundamental prophet
region, that is,
$$R_{\cal C}^{v,m} = \{(x,y) \; | \; x=V ({\bf X}),\; y=M ({\bf
X}); {\bf X} \in {\cal C}\} =  \{(x,y) \; | \; x \le y \le f_{\cal C}(x)\},$$ where
$f_{\cal C}$ is called the upper boundary function corresponding to ${\cal C}$.
Since it is only the latter set that gives a {\it complete} description of some
informational advantage, it is more fundamental and should be considered in its own
right.

In general, an {\it information set} characterizes the environment $\cal C$,
evaluated with the help of two particular levels of information. One could also
prioritise the information edge and say that the difference between two levels of information (e.g. minimal vs. sequential) is studied in a certain environment.
It is the second major aim of this article to illustrate a number of possible
applications of these ideas.

\section{Minimal versus maximum information}

In this section we systematically compare $u$ and $m$. That is, we are going to
derive corresponding information sets (called prophet regions since $m$ is involved)
in two standard random environments:
${\cal C}(I,n)$, the class of all sequences of independent, $[0,1]$-valued
random variables with horizon $n$; and ${\cal C}(G,n)$, the class of all
sequences of $[0,1]$-valued random variables with horizon $n$.

\begin{thm}\label{umi} (independent environment).
Let ${\bf X}=(X_1,\ldots,X_n) \in {\cal C}(I,n)$, $U({\bf X})=\max EX_i$ and
$M({\bf X})=E (\max X_i)$. Then the prophet region $\{(x,y) \; | \; x =U({\bf
X)}, y=M({\bf X}), {\bf X} \in {\cal C}(I,n)\}$ is precisely the set
$$R_{{\cal C}(I,n)}^{u, m} =\{ (x,y) | 0 \le x \le y \le f_n(x)=1-(1-x)^n;  \; 0 \le x \le 1 \} .$$
\end{thm}

\begin{thm} (general environment).
Let ${\bf X}=(X_1,\ldots,X_n) \in {\cal C}(G,n)$, $U({\bf X})=\max EX_i$, and
$M({\bf X})=E \max X_i$. Then the upper boundary function $h_n$ of the prophet
region $R_{{\cal C}(G,n)}^{u,m} = \{(x,y) \; | \; x =U({\bf X)}, y=M({\bf X}),
{\bf X} \in {\cal C}(G,n) \}$ is
$$
h_n(x) \:= \left\{ \begin{array}{lll} nx & \mbox{if} & 0 \le x < 1/n \\
1 & \mbox{if} & 1/n \le x \le 1 .
\\
\end{array}\right.
$$
\end{thm}

{\bf Proof of Theorem 1:} Without loss of generality let $x =EX_1\ge \max_{2 \le i
\le n} EX_i$. \citet[Lemma 2.2]{hi81} prove that ${\bf X}$ can be replaced by a
`dilated' vector ${\bf Y}$ of Bernoulli random variables $Y_1,\ldots,Y_n$ such that
$EX_i=EY_i$, $1\le i \le n$, and $M({\bf X}) \le M({\bf Y}) $. Replacing ${\bf Y}$
by a vector of iid Bernoulli random variables ${\bf Z}=(Z_1,\ldots,Z_n)$ such that
$EZ_i=x$, $1 \le i \le n$, does not improve the value to the gambler, i.e. $U({\bf
X})=U({\bf Y})=U({\bf Z})=x$, however, $M({\bf Y}) \le M({\bf Z}) =1-(1-x)^n$. Since
any ${\bf X} \in {\cal C}(I,n)$ can be replaced by a vector ${\bf Z}$ of iid
Bernoulli random variables without changing the value to the gambler, $f_n(x)$ is
the upper boundary function. Defining the independent random variables
$Z_1^{'},\ldots,Z_n^{'}$ by means of $P(Z_i^{'
}=\lambda+(1-\lambda)x)=x/(\lambda+x-\lambda x)=1-P(Z_i^{'}=0)$ and $0 \le \lambda
\le 1$ proves that all points between $(x,x)$ and $(x,1-(1-x)^n)$ also belong to the
region. $\diamondsuit$

Notice that for every fixed $x>0$, $\lim_{n\rightarrow \infty} f_n(x) \uparrow 1$
holds. Inspecting $f_n(x)/x$ and $f_n(x)-x$ immediately yields:

\begin{cor}
The prophet inequalities corresponding to ${\cal C}(I,n)$ are
$$M({\bf X}) / U({\bf X}) \le \lim_{x \rightarrow 0} f_n(x)/x =
n \;\;\mbox{and}\;\;\; M({\bf X}) - U({\bf X}) \le n^{-1/(n-1)}-n^{-n/(n-1)}.$$ In
the latter case, ${\bf Z}=(Z_1,\ldots,Z_n) \in {\cal C}(I,n)$ attains equality if
the $Z_i$ are iid Bernoulli random variables such that $U({\bf
Z})=EZ_i=P(Z_i=1)=1-n^{-1/(n-1)}$.
\end{cor}


\vspace{2ex}

{\bf Proof of Theorem 2:} Denote by $\bf e_i$ the $i$-th canonical unit vector.
First consider the random variable ${\bf Z} =(Z_1,\ldots,Z_n)$ having the
distribution
$$P({\bf Z}={\bf e}_1)= \ldots =P({\bf Z}={\bf e}_n) = 1/n.$$

A minimally informed person picks any of the random variables $Z_i$, which is 1 with
probability $1/n$ and obtains $U({\bf Z})=1/n$. Since there is always exactly one
$i$ such that $Z_i=1$, whereas all the other random variables are zero, $M({\bf
Z})=E(\max Z_i ) = \max_{1 \le i \le n} Z_i \equiv 1$.

To get a $U({\bf Z}) \ge 1/n$, let $P({\bf Z}={\bf e}_1)= x \ge 1/n$ and distribute
the remaining probability equally among the other canonical unit vectors, i.e.
$P({\bf Z}={\bf e}_2) = \ldots =P({\bf Z}={\bf e}_n) = (1-x)/(n-1) \le x$. Thus the
minimally informed gambler may always pick the first random variable, giving him
$U({\bf Z})=EZ_1=x$ and for the same reasons as before $E(\max Z_i) = 1$. Replacing
${\bf e}_i$ by $\lambda {\bf e}_i$ where $0 \le \lambda \le 1$ and $i=2,\ldots,n$
does not change the value to the gambler, but the value to the prophet decreases
towards $x$ if $\lambda \downarrow 0$.

Finally, let ${\bf X}=(X_1,\ldots,X_n)\in {\cal C}(G,n)$ and $U({\bf X})=x< 1/n$. On
the set $A_i =\{\omega | X_i (\omega) \ge \max_{j,j\neq i} X_j(\omega) \}$ replace
${\bf X}(\omega)=(X_1(\omega),\ldots,X_n(\omega))$ by

${\bf Y}(\omega)=(0,$ $\ldots,$ $0,$ $X_i(\omega),0,\ldots,0)$, $i=1,\ldots,n$. In
the case of equality choose any component (e.g. the first) where the maximum is
attained. Since $X_i(\omega)$ $\ge Y_i(\omega)$ we have $U({\bf X})=\max EX_i \ge
\max EY_i = U({\bf Y})=y,$ and since $\max [X_1(\omega),\ldots,X_n(\omega)] =
\max[Y_1(\omega),\ldots,Y_n(\omega) ]$, $M({\bf X}) = M({\bf Y})$.

By construction, at most one component of ${\bf Y}(\omega)$ is larger than zero.
Thus $\max_{1\le i \le n} Y_i(\omega) = \sum_{i=1}^n Y_i(\omega)$, and therefore

$$E(\max_{1 \le i \le n}Y_i)=E(\sum_{i=1}^n Y_i) = \sum_{i=1}^n EY_i \le
n \max_{1 \le i \le n} EY_i.$$

In the previous line equality is achieved if all expected values agree.
Defining the distribution of ${\bf Z} =(Z_1,\ldots,Z_n)$ via
$$P({\bf Z}={\bf e}_1)= \ldots =P({\bf Z}={\bf e}_n) = y \le x < 1/n \,\;\;\;
\mbox{and}\;\;\;\; P({\bf Z}={\bf 0})=1-ny$$ immediately yields $U({\bf Z})=y$ and
$M({\bf Z}) = ny.$ Since $y$ may assume any value in the interval $[0,1/n)$ we have
shown that $h_n(x)=nx$ is the upper boundary function if $x<1/n$. A similar
construction as before shows that all points between $(x,x)$ and $(x,nx)$ belong to
the prophet region. $\diamondsuit$

An immediate consequence of the last theorem is:
\begin{cor}
The prophet inequalities corresponding to ${\cal C}(G,n)$ are $M({\bf X}) / U({\bf
X}) \le n$ and $M({\bf X}) - U({\bf X}) \le 1-1/n$. In the latter case equality is
attained by $P({\bf Z}={\bf e}_i)=1/n$, $(i=1,\ldots,n)$ where ${\bf e}_i$ denotes
the $i$-th canonical unit vector.
\end{cor}

{\bf Remark.} Although we focus on the prophet, other comparisons, in particular
involving the statistician, would be interesting too. Comparing $u$ and $v$ for
example, reveals the difference between prior information on the one hand and
additional acquired information (sequential observations) on the other.

\section{Applying information sets}

In this section we restrict attention to classical prophet-statistician comparisons ($v$ vs. $m$).
However, the same kind of systematic analysis can be performed on any random
environment and observers with different levels of information. An example will
be given in the last section where we will compare $u$ an $m$.

\subsection{Some well-known results}

To illustrate how information sets may be used, we first collect a number of
well-known results. To this end we introduce further random environments: ${\cal
C}_{iid}$, the class of all sequences of iid, $[0,1]$-valued random variables;
${\cal C}_I$, the class of all sequences of independent, $[0,1]$-valued random
variables; ${\cal C}_G$, the class of all sequences of $[0,1]$-valued random
variables, and their corresponding counterparts with finite horizon i.e. ${\cal
C}^n_{iid}, {\cal C}^n_{I}={\cal C}(I,n)$ and ${\cal C}^n_{G}={\cal C}(G,n)$.

The $\beta$-discounted environment ${\cal C}_\beta$ is defined by $X_1=Y_1, X_2 =
\beta Y_2, X_3= \beta^2 Y_3,\ldots$, $(0 \le \beta \le 1)$ and ${\bf
Y}=(Y_1,Y_2,\ldots) \in {\cal C}_I$. Closely related are random variables
$X_1,\ldots, X_n$ with ``increasing bounds'', i.e. $a_i \le X_i \le b_i$ and
nondecreasing sequences $(a_i)$ and $(b_i)$. In both cases it suffices to study
$n=2$, i.e. $X_1 = \alpha Y_1, X_2=Y_2$ and $X_1=Y_1, X_2 = \beta Y_2$,
respectively, where $\alpha, \beta \in [0,1]$, and $(Y_1,Y_2) \in {\cal C}_I^2$.

The following table collects a number of well-known ``prophet'' results, i.e.
systematic comparisons of $v$ and $m$ (see \citet{hi83a, hi83, ke86, bo91}, and
\citet{sa98}):

$$
\begin{tabular}{|l|l|}
  \hline
  Random Environment & Upper boundary function \\\hline
${\cal C}_G$ & $f_{G}(x)=x-x \ln (x) $  \\
${\cal C}_G^n$ & $g_n(x)=nx-(n-1) \,x^{n/(n-1)} $ \\\hline

${\cal C}_{iid}$ & $x$  \\
${\cal C}_{iid}^n$ & $\phi_n(x)$ strictly increasing, strictly concave,
differentiable \\\hline

  ${\cal C}_I, {\cal C}^n_I$ & $f_I(x)=2x-x^2 $\\\hline

${\cal C}_\alpha^2$ & $2x-x^2$ if $x < \alpha$; and $x+(1-x) \alpha$ if $x \ge \alpha$ \\
\hline

 ${\cal C}_\beta, {\cal C}^n_\beta$ &   $f_\beta(x) =2x-x^2/\beta$ if $x <
 1-\sqrt{1-\beta}$, and \\
  & $f_\beta(x) =x+(1-x)(2(1-\sqrt{1-\beta})-\beta)$ if $x \ge 1-
\sqrt{1-\beta}$ \\\hline

\end{tabular}
$$

In general, the difficult part consists in finding an upper boundary function, yet it is easy
to show that all pairs $(x,y)$ with $x \le y < f_{\cal C}(x)$ belong
to some prophet region. 
Moreover, prophet inequalities follow straightforwardly from prophet regions. As
an example, look at $R_I$: Since $f_{I}(x) / x = 2-x \le 2$ and $f_{I}(x)- x =x
(1-x) \le 1/4$, we have $M({\bf X})/V({\bf X}) \le 2$ and $M({\bf X})-V({\bf X})
\le 1/4$ for all ${\bf X} \in {\cal C}_I$. The same kind of argument yields
$M({\bf X}) < V({\bf X})(1- \ln V({\bf X}))$ and $M({\bf X})-V({\bf X}) < 1/e$
for all ${\bf X} \in {\cal C}_G$.


\subsection{Graphical comparisions}

What can be learned from this upon comparing two gamblers with different information
levels? For every fixed horizon $n$, we have $R_{iid}^n \subseteq R_I^n \subseteq
R_G^n$. It also turns out that $R_{iid}^n \subset R_{iid}^m$ and $R_{G}^n \subset
R_{G}^m$ whenever $n < m$. Thus the longer the horizon or the more general the
environment, the better the outcome for the prophet (or the better informed person
in general). On the other hand, restrictions of any kind, in particular the range of
the random variables makes the corresponding prophet (or information) region
smaller. For example, $R_\alpha^2$ and $R_\beta$ must be subsets of $R_I$.

The following illustration combines results achieved so far.

\pagebreak
\begin{figure}[h]

\includegraphics{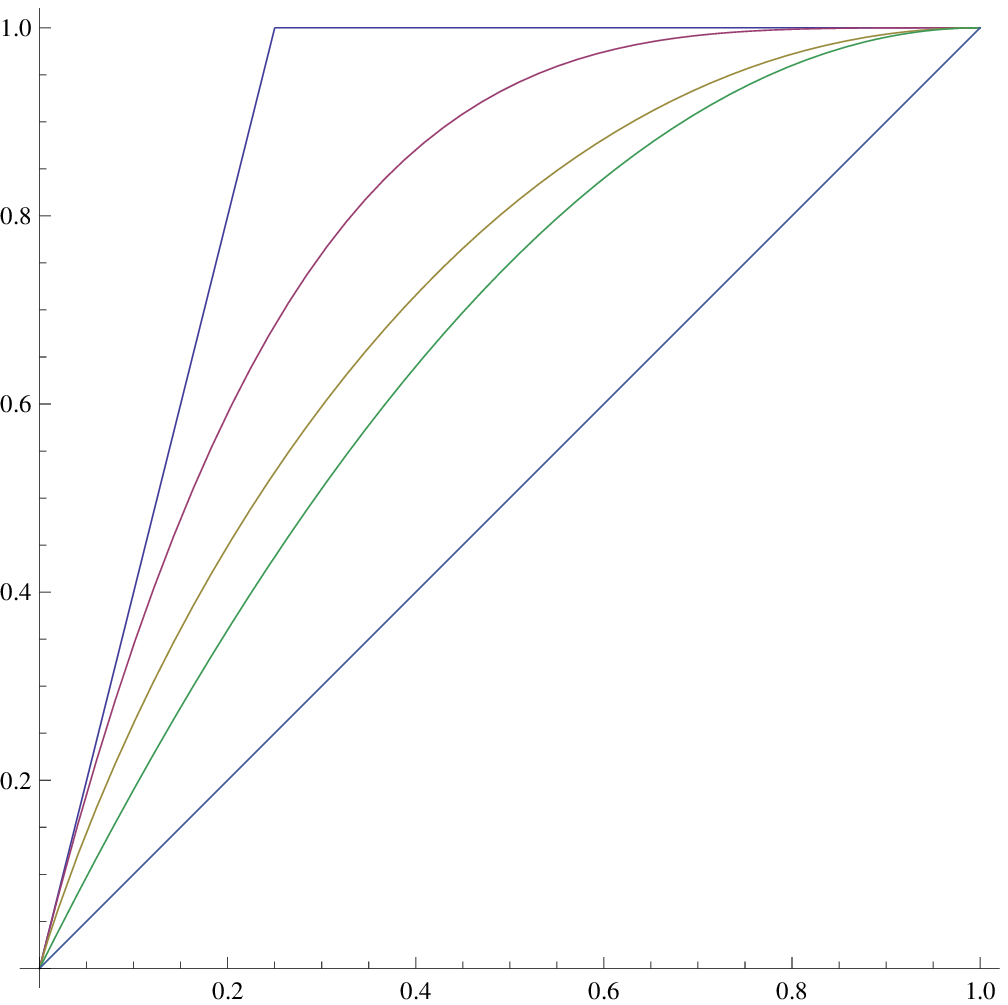}

\end{figure}

{\bf Illustration 1.}\label{family} From above: The functions
$h_4(x), f_4(x), g_4(x), f_I(x)$, and $x$.

Note that $h_4$ and $f_4$ stem from comparisons of $u$ and $m$, whereas $g_4$ and
$f_I$ are the result of comparisons of $v$ and $m$ in the general and the
independent environments. Since for any environment $R_{\cal C}^{v,m} \subseteq
R_{\cal C}^{u,m}$, we must have $g_4 \le h_4$ and $f_I \le f_4$. In the case $n=2$
the functions $f_I$ and $f_2$ agree. This is no coincidence since $X_1 \equiv x$ and
$x=P(X_2=1)=1-P(X_2=0)$ is the (standard) worst case scenario for the statistician,
and $x=P(X_i=1)=1-P(X_i=0)$ $(i=1,2)$ is the worst case scenario for the minimally
informed gambler considered above. In both scenarios their values agree (e.g. they
may both choose the second random variable) giving the prophet a maximum advantage
of $x(1-x)$.

\subsection{The overall information difference}

The diagonal `$y=x$' collects all situations where the information edge of a better
informed person does not result in a larger payoff. Thus, a degenerated prophet
region indicates that given a stochastic environment the information lead of the
prophet never pays off. Yet, the further some upper boundary function is away from
the identical function, the larger the better informed gambler's overall advantage.
A natural measure of this advantage is the area between these functions, i.e. the
integral
$$
\int_{0}^1 (f_{\cal C}(x)-x) \; dx .
$$
Given ${\cal C}_I$, the prophet's advantage is $q_I=\int_{0}^1 x \,(1-x) \; dx =
1/6$. In the discounted environment, after some algebra, we obtain
$$
q(\beta)
 = \frac{1}{6} - \frac{(1-\beta) \; (1-\sqrt{1-\beta})}{3 \beta} .
$$
Note that $q(1)=a_I=1/6$, and l`Hopital's rule gives $\lim_{\beta \downarrow 0}
q(\beta) =0$. Moreover, $q(\beta)$ is a convex function.

In the ``increasing bounds'' environment, after a little algebra, we obtain
$$
{\tilde q}(\alpha)  = \int_0^{\alpha} (2x-x^2 ) \;dx +\int_\alpha^1 (\alpha-\alpha x
+x)\; dx -\frac{1}{2} = \alpha (\alpha^2 /3-\alpha  + 1) /2.
$$

Note that ${\tilde q}(1)=q_I=1/6$, and l`Hopital's rule gives $\lim_{\alpha
\downarrow 0} {\tilde q}(\alpha) =0$. Moreover, ${\tilde q}(\alpha)$ is a concave
function.

\begin{figure}[h]

\includegraphics{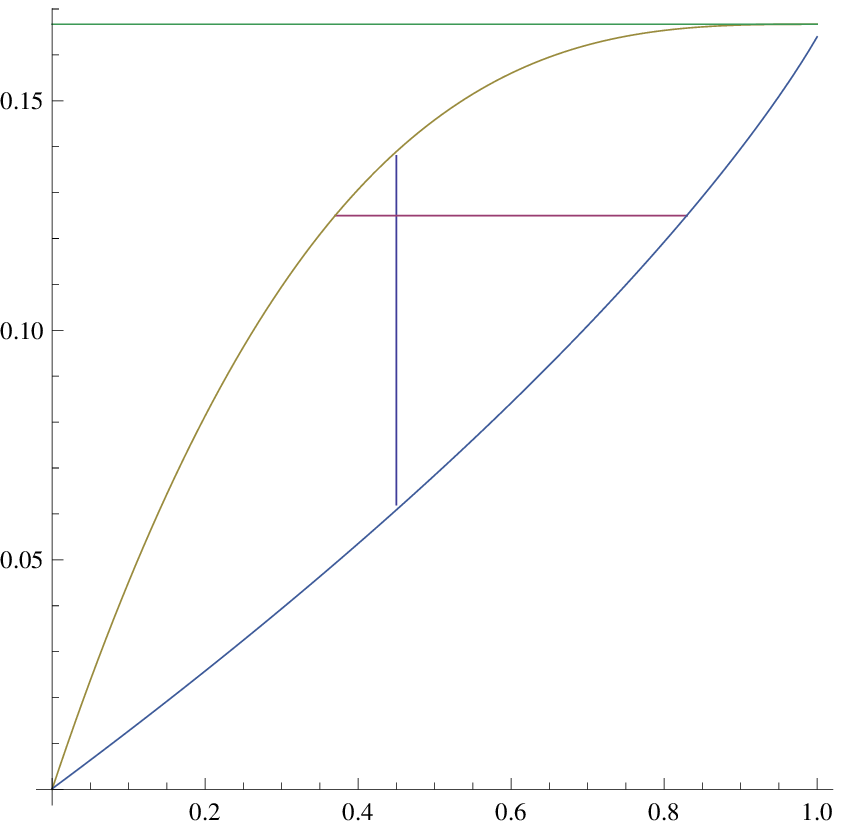}

\end{figure}

{\bf Illustration 2.} $\alpha$ and $\beta$ are shown on the $x$-axis. The functions
on the unit interval from the top down are the constant $q_I=1/6$, ${\tilde
q}(\alpha)$ and $q(\beta)$. The vertical and the horizontal lines will be explained
in Section \ref{environ}.

However, given ${\cal C} (G,n)$, $q_I$ is augmented to
$$q_G(n)=\int_{0}^1 (n-1) \,(x-x^{n/(n-1)}) \; dx
=\frac{n-1}{2 \,(2n-1)} .$$ In particular, $q_G(2)=1/6$, $q_G(3)=1/5$,
$q_G(4)=3/14$, $q_G(5)=2/9$, and $q_G(6)=5/22$. Moreover, $q_G(n)$ is strictly
increasing in $n$ with limit $1/4$, and
$$q_G(\infty)= -\int_{0}^1 x \ln(x) \; dx
=- \left[ x^2  \left( \frac{\ln(x)}{2} - \frac{1}{4} \right)  \right]_0^1 =
\frac{1}{4}.$$

\subsection{Inverse problems}\label{inv}

Given a stochastic environment ${\cal C}$, and according to the above derivation,
the standard interpretation of a prophet inequality, is to look for a value to the
statistician $x_0=V({\bf X})$, such that the difference $f_{\cal C}(x)-x$ is
maximized. In the same vein one may look for a value $y_0$ on the $y$-axis, where
the difference between the upper boundary and the identical function is at its
greatest point. In the independent case this amounts to inverting $f_I(x) = 2x-x^2$,
which yields $f_I^{-1}(y)=1-\sqrt{1-y}$. Maximizing $y-(1-\sqrt{1-y})$ gives $1/4$,
which is obtained for $y_0=3/4$.

Why do both perspectives agree with respect to the maximum difference? The
reason is that the statement $M({\bf X})-V({\bf X}) \le 1/4$ holds for all ${\bf
X} \in {\cal C}_I$, and thus is a property of the stochastic environment (and
the two levels of information considered). The pair $(1/2, 3/4) \in R_I$ is a
point in two-dimensional space, attained by certain extremal sequences ${\bf
X}^*$. Thus, no matter how we choose to look at some region $R_{\cal C}$, the
corresponding prophet inequalities must hold.

However, the analytic considerations involving the inverse of the upper boundary
function may be quite different. 
In the discounted case, $f_\beta^{-1}(y)=\beta(1-\sqrt{1-y/\beta})$ if $y \le
g(1-\sqrt{1-\beta})=3-2/\beta -2\sqrt{1-\beta}+2\sqrt{1-\beta}/\beta= y(\beta)$.
Otherwise, it is easily seen that $f_\beta^{-1}$ is a linear, strictly
decreasing function of $y$, and $f_\beta^{-1}(1)=1$. The maximum of the function
$y-\beta(1-\sqrt{1-y/\beta})$ occurs at the point $y=3\beta/4$ and is $\beta/4$.
Notice that
$$y^{'}(\beta) = \frac{\beta^2 +\beta-2+2\sqrt{1-\beta}}{\beta^2 \sqrt{1-\beta}}
\ge  \lim_{\beta \downarrow 0} \frac{\beta^2 +\beta-2+2\sqrt{1-\beta}}{\beta^2
\sqrt{1-\beta}} = \frac{3}{4}.$$ Thus, $3\beta/4 < y(\beta)$ for all $\beta > 0$.
Due to continuity of $f_\beta^{-1}$, this yields $\beta/4$ as the overall maximum of
the difference, always occuring at $y=3\beta /4$. Traditionally, one would have said
that the maximum difference of $\beta/4$ occurs at $x=\beta/2$.

Given ${\cal C}_G$, one has to invert $f_G(x)=x-x\ln(x)$ in the unit interval. Using
the theorem of the derivative of the inverse function one may check that
$\exp(1+W_{-1}(-y/e))$ is the inverse, where $W_{-1}(y)$ is the lower real branch of
the Lambert $W$ function \citep[see][331]{co96}. Thus $(y-\exp(1+W_{-1}(-y/e)))' =
1+1/(1+W_{-1}(-y/e))$ and $W_{-1}(-2/e^2)=-2$ immediately yield that the maximum
occurs at $y=2/e$ and equals $1/e$. Traditionally, it's the same difference occuring
at $x=1/e$.

\subsection{Comparing stochastic environments}\label{environ}

Switching stochastic environments amounts to a systematic comparison of the
associated regions. In particular, if ${\cal A}$ is less general than ${\cal
B}$, we have $R_{\cal A} \subseteq R_{\cal B}$. Obviously, it suffices to
consider the upper boundary functions $f_{\cal A}, f_{\cal B}$ of the two
environments involved. Traditionally, one would only determine $\sup_x ( f_{\cal
B} (x) - f_{\cal A} (x) )$. However, the inverse problem $\sup_y ( f_{\cal
A}^{-1} (y) - f_{\cal B}^{-1} (y) )$, and the area $\int_0^1 (f_{\cal B}(x)-
f_{\cal A}(x)) \; dx$ are also natural measures of discrepancy.

To illustrate the above, let us compare ${\cal C}_I$ and ${\cal C}_G$:

First, $\int_0^1 (f_G(x)-f_I(x)) \; dx=3/4-2/3=1/12$.

Second, maximizing $d(x)=f_G(x)-f_I(x)= x^2-x-x \ln x$ leads to $d^{\,'}(x)=0
\Leftrightarrow 2x-\ln x = 2$, which has the explicit solution $x_0=-W_0(-2/e^2)/2
\approx 0,406376/2$, where $W_0$ is the principal (upper) real branch of the Lambert
$W$ function \citep[see][331]{co96}. The point $(x_0,d(x_0))\approx(0.2,0.162)$ may
be interpreted as follows: For every value $x$ to the statistician, $f_I(x)$ is the
best a prophet can obtain in the independent environment ${\cal C}_I$, and he can
get arbitrary close to $f_G(x)$ if he is confronted with the general environment
${\cal C}_G$. Given $x$, the difference $f_G(x)-f_I(x)$ reflects the additional gain
(almost) obtainable to the prophet when moving from ${\cal C}_I$ to ${\cal C}_G$,
i.e. from the restricted to the more general situation. The additional sequences of
random variables provide him with an additional reward of $d(x)=x(x-\ln x -1)$,
which is maximized if $x=-W_0(-2/e^2)/2$, yielding $0.162$ as the additional payoff.

Third, starting with the prophet, the difference to be considered is $\delta(y)
=f_I^{-1}(y)-f_G^{-1}(y)=1-\sqrt{1-y}-\exp (1+W_{-1}(-y/e))$. Thus, conditional on
$y$, the statistician may (almost) lose this amount when the stochastic environment
switches from independent to arbitrary sequences of random variables. Determining
the value $y_0$ where $\delta(y)$ is at its greatest, means looking for a
constellation where the loss occuring to the statistician is the most pronounced
when moving from ${\cal C}_I$ to ${\cal C}_G$. Now $\delta^{\, '}(y)=0$ is
equivalent to finding the unique root of the equation
$$
-2\sqrt{1-y} = 1+W_{-1}(-y/e) .
$$
As a function of $y$, both the left hand side (L) and the right hand side (R) of the
equation are twice differentiable. On the unit interval $L(y)$ is convex, strictly
increasing, $L(0)=-2$, and $L(1)=0$. $R(y)$ is concave, strictly increasing,
$\lim_{y \downarrow 0} R(y) = -\infty$, and $R(0)=0.$ Numerically, this yields the
solution $(y_0,\delta(y_0)) \approx (0.70,0.119)$. Thus, in the worst case, the
statistician loses about $0.119$, which is considerably less than the prophet can
hope to obtain when the environment extends from ${\cal C}_I$ to ${\cal C}_G$.

The next illustration summarizes these results:

\begin{figure}[h]

\includegraphics{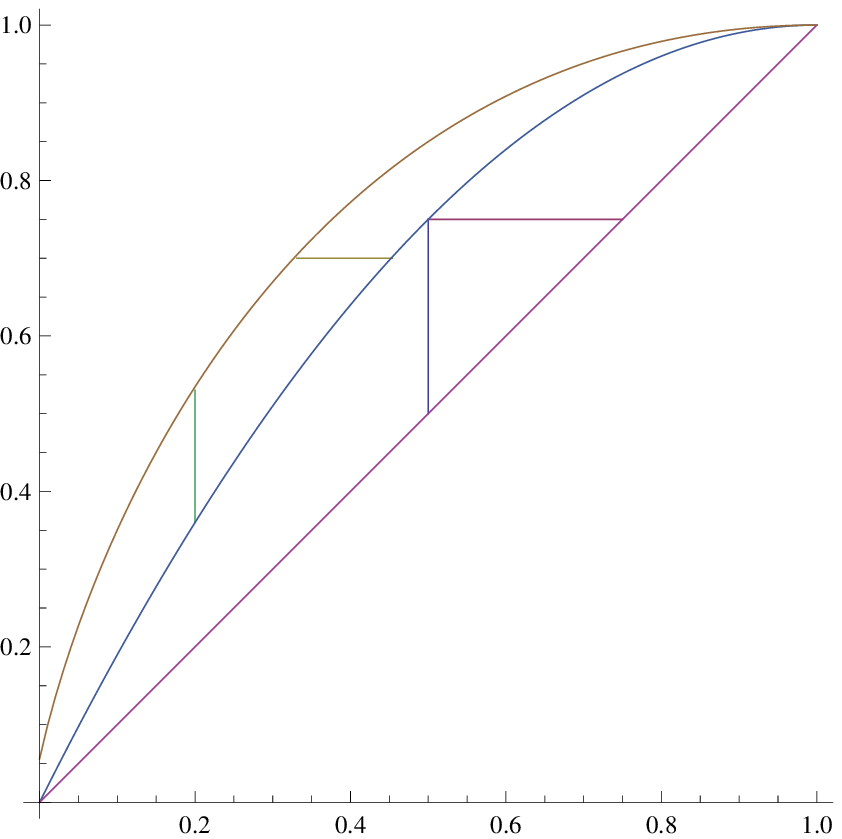}

\end{figure}

{\bf Illustration 3.} From above: The functions $f_G, f_I$, and $x$ on
the unit interval. The small vertical line to the left illustrates the position of
the maximum of the function $d(x)$, the small horizontal line illustrates the
maximum of the function $\delta(y)$, see Section \ref{environ}. The other lines
indicate the position of the maximum difference between the statistician and the
prophet in the independent environment, see the second paragraph of Section
\ref{inv}.

A different kind of analysis may be explicated using the regions ${\cal C}_\alpha^2$
and ${\cal C}_\beta^2$:
Illustration 2 points out that restricting the range of the second random variable
($\beta$-discounting), always produces a smaller region than restricting the range
of the first random variable by the same amount $\alpha$.
The largest difference between the size of the regions occurs if $\alpha=\beta
\approx 0.45$ and is approximately $0.077$. On the other hand suppose that the areas
of $R_\alpha$ and $R_\beta$ agree. This is tantamount to fixing a point on the
$y$-axis. In this case the largest difference between the parameter values occurs if
the area covered by each of the regions is about $1/8$. There $\alpha \approx 0.38$
and $\beta \approx 0.83$, thus the largest difference between the parameter values
is approximately $0.452$.

Of course, analyses along the same lines can be carried out for other regions,
e.g., $R_I$ and $R_{iid}^n$, $R_{iid}^n$ and $R_{iid}^{n+1}$, $R_{G}^n$ and
$R_{G}^{n+1}$, or $R_G^n$ and $R_G$.

\subsection{Typical differences and ratios}

Classical prophet inequalities are `worst case' scenarios. They refer to the
maximum advantage of the prophet over the statistician. Additionally, it is
straightforward to ask for a `typical' advantage, in particular a `typical'
difference or ratio. To do so, one would have to define a probability measure on
some environment $\cal
 C$. Since the classes of random variables considered are rather large, it is by no
means clear how to do so in a natural way. However, starting with a stochastic
environment {\it and} two distinguished levels of information, it is natural to
consider uniform measure on the corresponding prophet region $R_{\cal C}$.

Given the independent environment, the size of $R_I$ is $1/6$. Thus, we obtain
as the typical difference between $M({\bf X})$ and $V({\bf X})$
$$
6 \int_{R_I} (y-x) \;d(x,y) = 6 \int_0^1  (x^4/2-x^3+x^2/2) \; dx  = 1/10 ,
$$
instead of $1/4$ in the worst case. Moreover, the typical ratio is
$$
6 \int_{R_I} y/x \;\;d(x,y) = 6 \int_0^1 (x^3/2 -2x^2+ 3x/2) \; dx =
5/4 < 2.
$$

Given ${\cal C}_G$, $R_G$ covers an area of $1/4$, giving the
following typical difference and ratio:
$$
4 \int_{R_G} (y-x) \;\;d(x,y) = 4 \int_0^1 (x^2 \ln^2(x) /2) \; dx = 4/27 < 1/e
$$
and
$$
4 \int_{R_G} y/x \;\;d(x,y) = 4 \int_0^1 (x \ln^2(x) /2 -x \ln(x)) \; dx = 3/2.
$$
The last equation is particularly interesting because there is no upper bound in
the corresponding worst case scenario. Notice in the other examples that the
typical results are considerably smaller than the constants in the corresponding
worst cases.

Moreover, one may ask about the probability that a typical difference or ratio
exceeds a certain bound. The ratio $y/x =c \Leftrightarrow y = cx$ is a straight
line through the origin, so, given ${\cal C}_I$, the question amounts to
calculating
$$
P(\mbox{M({\bf X})/V({\bf X}))} \ge c) = 
\int_0^t (2x-x^2 -cx) \; dx 
/ \; \int_0^1 (2x-x^2 -x) \; dx =(2-c)^3 ,
$$ where $t=2-c \ge 0$ is determined by
the equation $cx=y=2x-x^2$, and $1 \le c \le 2$. Given ${\cal C}_G$, we obtain
$$
P(\mbox{M({\bf X})/V({\bf X}))} \ge c) = 
\int_0^t (x-x \ln x -c x) \; dx 
/ \; \int_0^1 (x-x \ln(x) -x) \; dx = e^{2(1-c)} ,
$$
where $t$ is determined by the equation $cx=x-x \ln x \Leftrightarrow t =
\exp(1-c)$, and $c \ge 1$.

In the case of the difference we are interested in the probability that it
exceeds a certain bound $d \ge 0$. Again, consider ${\cal C}_I$ first. Since
$y-x =d \Leftrightarrow y = x+d$, we have to calculate
$$
P(M({\bf X})-V({\bf X}) \ge d) = \frac{\int_s^t (2x-x^2 -(x+d)) \; dx }{
 \int_0^1 (2x-x^2 -x) \; dx} = \sqrt{1-4d} \cdot (1-4d) .
$$
Here, $0 \le d \le 1/4$, and $s$ and $t$ are determined by the roots of the
equation $x+d=2x-x^2$ in the unit interval, that is $s=1/2-\sqrt{1-4d}/2$ and
$t=1/2+\sqrt{1-4d}/2$.

Finally, given ${\cal C}_G$, we obtain with $0 \le d \le 1/e$
$$
P(D \ge d)=P(M({\bf X})-V({\bf X}) \ge d)  = 
-\int_{s}^{t} (x \ln x +d) \;dx\;\; 
/\; a_G(\infty),
$$
where $s$ and $t$ are determined by the roots of the equation $d=-x \ln x$ in
the unit interval. Some algebra is needed to get $s=\exp (W_{-1}(-d))$ and
$t=\exp (W_{0}(-d))$. The subsequent integration results in
$$
P(D \ge d)  =  d^2\left( \frac{1+4W_0(-d)-2
\ln\left(\frac{-d}{W_0(-d)}\right)}{W_0^2 (-d)}
   -  \frac{1+4W_{-1}(-d)-2 \ln\left(\frac{-d}{W_{-1}(-d)}\right)}{W_{-1}^2
   (-d)} \right) .
$$

\section{A systematic study}

In the following we are going to apply the `program' outlined in the last section to $u$ and $m$, using the independent and the general stochastic environments:

\vspace{2ex}
{\bf The overall information difference}. Let us first compute the areas of
$R^{u,m}_{{\cal C}(I,n)}$,
$$
\int_0^1 (f_n(x)-x) \; dx = \frac{1}{2}-\frac{1}{n+1} = \frac{n-1}{2(n+1)}
\rightarrow 1/2,$$ and $R^{u,m}_{{\cal C}(G,n)}$;
$$\int_0^1 (h_n(x)-x) \;
dx = \frac{1}{2n}-\frac{1}{2n^2}+\frac{1}{2}-\frac{1}{n}+\frac{1}{2n^2} =
\frac{1}{2}-\frac{1}{2n}=\frac{n-1}{2n} \rightarrow 1/2 .
$$
Thus, their overall information distance is the size of the set $R^{u,m}_{{\cal
C}(G,n)} \backslash R^{u,m}_{{\cal C}(I,n)}$,
$$
\int_0^1 (h_n(x)-f_n(x)) \; dx = \frac{1}{n+1}-\frac{1}{2n} =
\frac{n-1}{2n(n+1)} \rightarrow 0 \mbox{ if } n \rightarrow \infty.
$$

{\bf Inverse Problems}. In the independent case, $f^{-1}_n(y) = 1- \sqrt[n]{1-y}$ is
the inverse function. The maximum of $y-(1- \sqrt[n]{1-y})$ is attained for $y_0 =
1-n^{-n/(n-1)}$ and equals $n^{-1/(n-1)}-n^{-n/(n-1)}$. In the general case, the
inverse function is $h^{-1}_n(y) = y/n$. Thus, the maximum of $y-y/n$ is attained at
$y_0=1$, giving a maximum difference of $1-1/n$.

{\bf Comparing the independent and the general environments}. Here, one has to
maximize $d(x) = (1-x)^n + nx -1$. Since $d^{'}(x) = n(1-(1-x)^{n-1}) >0$ if $x \le
1/n$, the maximum occurs at $x_0 = 1/n$ and yields a difference of $(1-1/n)^n$,
converging to $1/e$ if $n\rightarrow \infty$. The inverse functions lead to a
difference of $\delta(y)=1-(1-y)^{1/n}-y/n$. Thus, $\delta^{'}(y) =
-(1-(1-y)^{(n-1)/-n})/n > 0$ if $y >0$, and the maximum occurs at $y_0 = 1$,
yielding a difference of $1-1/n$.

{\bf Typical differences and ratios}. For ${\cal C}_I^n$ we calculate
$$
\int_{R_{{\cal C}(I,n)}^{u,m}} (y-x) \;\;d(x,y) = \frac{1}{6}-\frac{1}{n+2}+
\frac{1}{2(2n+1)}=\frac{(n-1)^2}{3(n+2)(2n+1)}
$$
and
$$
\int_{R_{{\cal C}(I,n)}^{u,m}} y/x \;\;d(x,y) =-\frac{1}{4}+\sum_{i=1}^n (1/i) -
\frac{1}{2} \sum_{i=1}^{2n} (1/i).
$$
Thus, the typical difference $d_I$ and ratio $r_I$ in the independent situation
are
$$d_I=\frac{2(n-1)(n+1)}{3(n+2)(2n+1)} \rightarrow 1/3,
\;\;r_I=\frac{(n+1)(2\sum_{i=1}^n (1/i) -  \sum_{i=1}^{2n} (1/i)-1/2)}{(n-1)}
\rightarrow \infty .
$$

For ${\cal C}_G^n$, analogous integrations yield
$$
\int_{R_{{\cal C}(G,n)}^{u,m}} (y-x) \;\;d(x,y) = \int_0^{1/n} \int_{x}^{nx}
(y-x) \; dy \; dx +\int_{1/n}^1 \int_x^1 (y-x) \; dy \; dx =
\frac{(n-1)^2}{6n^2}
$$
and
$$
\int_{R_{{\cal C}(G,n)}^{u,m}} y/x \;\;d(x,y) = \int_0^{1/n} \int_{x}^{nx} y/x
\; dy \; dx +\int_{1/n}^1 \int_x^1 y/x \; dy \; dx = \frac{\ln n}{2}.
$$
Thus, the typical difference $d_G$ and ratio $r_G$ in the general environment
are
$$d_G =\frac{n-1}{3n} \rightarrow 1/3 \;,
\;\;r_G =\frac{n \ln n}{n-1} \rightarrow \infty .
$$

{\bf Probabilities that a typical difference or ratio exceeds a certain bound}. For
$1 \le c \le n$ this amounts to calculating
\begin{eqnarray*}
P(r_I \ge c) &=& \frac{2(n+1)}{n-1} \int_0^t (1-(1-x)^n -cx) \;dx
\\
&=& \frac{2}{n-1} \left( (1-t)^{n+1} + (n+1) t \left(1-\frac{ct}{2}\right) -1
\right) ,
\end{eqnarray*}
where $t$ is the unique root of the equation $cx=1-(1-x)^n$ in the unit
interval, and  
$$
P(r_G \ge c) =\frac{2n}{n-1}
\left( \int_0^{1/n} (nx-cx) dx + \int_{1/n}^{1/c} (1-cx) dx \right)
  = \frac{n-c}{c(n-1)} .
$$
Notice that $\lim_{n \rightarrow \infty} (n-c)/(c(n-1)) = 1/c$.

In the case of the difference, given ${\cal C}_I^n$, and thus $0 \le d \le
n^{-1/(n-1)}-n^{-n/(n-1)}$, we calculate
$$
P(d_I \ge d) = 
\frac{2(n+1)}{n-1} \int_s^t (1-(1-x)^n -(x+d)) \; dx  , 
$$
where the values of $s$ and $t$ $(s<t)$ are determined by the roots of the equation
$x+d=1-(1-x)^n$ in the unit interval. Again, in general, $s$ and $t$ cannot be given
explicitly. Finally, given ${\cal C}_G^n$, we obtain with $0 \le d \le 1-1/n$
\begin{eqnarray*}
P(d_G \ge d) &=& 
\frac{2n}{n-1} \left( \int_{a/(n-1)}^{1/n} (nx-x-d) \; dx + \int_{1/n}^{1} (1-x-d) \; dx \right) 
  \\
&=& 1- \left( \frac{dn}{n-1} \right) \left( 2-d-\frac{d}{n-1}\right) \rightarrow
(1-d)^2 \;\mbox{if} \; n \rightarrow \infty.
\end{eqnarray*}

In both cases the prophet regions ${R_{{\cal C}(I,n)}^{u,m}}$ and ${R_{{\cal
C}(G,n)}^{u,m}}$ converge towards the upper triangle $T=\{(x,y) | 0 \le x \le y
\le 1\}$ in the unit square. Thus, in the limit, the typical ratios and
differences agree and can be computed directly via $T$, yielding the
probabilities $1/c$ and $(1-d)^2$.

\end{document}